\documentstyle[12pt,a4]{article}
\title{Wedge removability of metrically thin sets and application to
  the CR-meromorphic extension}
\author{Tien-Cuong Dinh and Frédéric Sarkis }
\date{ \ }
\def\fdem {{$\hfill\Box$}}
\newtheorem{defin}{Definition}

\newtheorem{prop}{Proposition}
\newtheorem{lemme}{Lemma}
\newtheorem{corol}{Corollary}
\newenvironment{demo}{{\em Proof.}}{\fdem}
\def\k#1{\kern#1em}
\def\Ib#1{{I\kern-.25em#1}}
\def\Ibb#1{{I\kern-.23em#1}}

\def\vcg{\vrule  width.02em height1.4ex depth-.05ex}

\def\CC{{\rm{\k{.24}\vcg\k{-.26}C}}}

\def\LL{{\rm\Ib L}}

\def\NN{{\rm\Ibb N}}

\def\RR{{\rm\Ib R}}

\def\Im{{\rm Im}}
\def\Re{{\rm Re}}

\begin{document}
\maketitle
{\center{\bf Abstract.}}\par\noindent
We give a wedge removability theorem for metrically thin sets of two
codimensional Hausdorff null measure. Following  \cite{fredy},
this removability theorem combined with the wedge removability theorem
of \cite{merker2} for closed subsets of two codimensional manifolds,
gives a CR-meromorphic extension theorem in the greater
codimensional case.
\medskip
\section{Introduction}
\noindent {\it 1.1 CR submanifolds - }
For a smooth submanifold $M$ of an open subset of
$\CC^n$, let $T_p(M)$ be the real tangent
space of $M$ at $p \in M$. In general, $T_p(M)$ is not invariant under
the complex structure map $J$ for $T_p(\CC^n)$. Therefore, we give
special designation for the largest $J$-invariant subspace of $T_p(M)$.
For a point $p \in M$, the {\it complex tangent space} of $M$ at $p$ is the
vector space $$H_p(M)=T_p(M) \cap J\{T_p(M)\}.$$
The {\it totally real part} of the tangent space of $M$ is the quotient
space $$X_p(M)=T_p(M)/ H_p(M).$$
The complexifications of $T_p(M)$, $H_p(M)$ and $X_p(M)$ are denoted
by $T_p(M) \otimes \CC$, $H_p(M) \otimes \CC$ and $X_p(M) \otimes
\CC$. The complex structure map $J$ on $T_p(\RR^{2n})$ restricts to a
complex structure map on $H_p(M) \otimes \CC$ because $H_p(M)$ is
$J$-invariant. So, $H_p(M) \otimes \CC$ is the direct sum of the $+i$
and $-i$ eigenspace of $J$ which are denoted by $H_p^{1,0}(M)$ and
$H_p^{0,1}(M)$.
A submanifold $M$ in an open subset $U$ of $\CC^n$ is called an {\it
embedded
CR manifold or a CR submanifold} of $U$ if $\dim_\CC H_p(M)$ is independent
of $p \in M$.
A CR submanifold $M$ is called {\it generic} if {\it the CR dimension}
$\dim_\CC H_p(M)$ is minimal.\par
\medskip
\noindent{\it 1.2 CR and CR-meromorphic functions -}
Let $M$ be a CR submanifold of an open subset of $\CC^n$. A function $f: M
\rightarrow
\CC$ is called a {\it CR function }if $\bar \LL f = 0$ (in the current
sense on $M$) for any $C^1$ vector field $\LL$ of $H^{1,0}(M).$
A closed subset $N$ of a manifold $X$ is called a {\it scarred manifold} of
dimension $p$ if there exists a closed subset $\tau \subset N$
($\tau$ is called the {\it scar set}) of $p$-dimensional Hausdorff
${\cal H}^p$ zero measure such
that $N \backslash \tau$ is an oriented $p$ dimensional $C^1$ submanifold
of $X \backslash \tau$ of locally finite ${\cal H}^p$ volume in $X$ and
closed
in the current sense ($d[N]=0$). If the regular part of $N$ is a
maximally complex or a CR submanifold, we will say that $N$ is
a scarred maximally complex or CR submanifold.
\begin{defin}\rm (\cite{fredy})
Let $M$ be a CR manifold of dimension $p$.
A  $C^1$ CR map $f$ defined on an open and dense subset of $M$
 and having values in $P_1(\CC)$ will be called a {\it CR-meromorphic
   map}  if the closure $\Gamma_f$ of its graph in
$M \times P_1(\CC)$ is a $C^1$ scarred CR manifold of the
same CR dimension than $M$.
\end{defin}
This notion was first introduced by Harvey and Lawson
\cite{harvey2} in the case $M$ was maximally complex.
In that case, the meromorphic extension of CR-meromorphic maps  is obtained
by solving the boundary problem for their graph (see \cite{harvey2}
\cite{fredy}). In the case of non maximally complex CR manifold,
the solving of the boundary problem seems not to apply.
According to \cite{fredy}, in the case that a CR-meromorphic map
take its values in a bounded domain of $\CC$, it defines a
CR current on $M$. Thus, from the representation theorem of Baouendi and
Rothschild \cite{baoundi} we are reduced to a problem
of extension of smooth CR functions defined on $M$.
Let $f$ be a CR-meromorphic function defined on $M$ and let
$x \in M$ be a point such that there exists a point $y \in P_1(\CC)$
such that $x \times y \not \in \Gamma_f$.  Let $\phi_y: P_1(\CC)\backslash
y \rightarrow \CC$ be the projective
chart where $y$ is the infinite point and $\omega$ be a small enough
neighborhood of $x$ in $M$ such that $\phi_y \circ f$ is bounded.
Then $\phi_y f$ defines a CR current on $\omega$.
A point $p \in M$ is said an {\it indeterminacy
point} of $f$ if $\Gamma_f \supset \{p\} \times P_1(\CC)$.
The {\it indeterminacy set} $K$ of $f$ is the set of all indeterminacy
points of $f$. $K$ is an obstruction to the extension. Indeed,
if $x \in K$, the reduction to the case of CR currents cannot apply.
\par
\medskip
\noindent {\it 1.3 The Levi form - }
The local extension of CR functions or CR currents arise under
convexity assumption of $M$. The more general notion of convexity
where the extension arise is the notion of local and global
minimality. However, as a first step lets recall
to the notion of convexity given by the Levi form.\par\noindent
Suppose $M = \{ \zeta \in \CC^n; \rho_1(\zeta)=...=\rho_d(\zeta)=0\}$
is a smooth CR submanifold of $\CC^n$, with $1 \leq d \leq n$. Let $p$
be a point in $M$ and suppose $\{ \bigtriangledown
\rho_1(p),...,\bigtriangledown \rho_d(p)\}$ is an orthonormal basis
for the normal space $N_p(M)$ of $M$. Then the {\it extrinsic Levi form} is
given by $$ \widetilde {\cal L}_p(W)=-\sum_{l=1}^{d}
\left(
\sum_{j,k=1}^n { \partial^2 \rho_l(p)\over \partial \zeta_j
  \partial \bar \zeta_k } w_j \bar w_k \right) \bigtriangledown \rho_l(p)$$
for $W=\sum_{k=1}^n w_k(\partial / \partial \zeta_k ) \in H_p^{1,0}(M)$.\par
\medskip
\noindent{\it 1.4 Minimality - }
One of the characteristic properties of CR manifold is that
$H(M)$ is involutive. So, we will say that a curve
$\gamma:[0,1] \rightarrow M$ is a {\it CR curve} if for all $t \in
[0,1]$, the tangent to $\gamma$ at the point $\gamma(t)$ is in
$H_{\gamma(t)}(M)$. We will say that two points are in the same {\it
  CR orbit of M} if they can be reached by a finite number of CR curves of
$M$. Let $\{U_i\}_{i \in I}$ be a basis of neighborhoods of $p$ in $M$. For each
$i\in I$, we can define the orbit of $p$ in $U_i$. The inductive
limit of those orbit is well defined and does not depends of the
family $\{U_i\}$ and is call the {\it local CR orbit} of $p$.
A CR submanifold $M$ is said minimal at $p$ is the local CR orbit of
$p$ is a open neighborhood of $p$ in $M$.
\medskip
\par\noindent
{\it 1.5 Wedge removable sets - }
Let $M$ be a generic CR submanifold of an open subset of $\CC^n$ (i.e. $T_pM+JT_pM=T_p\CC^n$).
By a {\it wedge of edge $M$ at $p$}, we mean an open set in $\CC^n$ of the
form $${\cal W}=\{ z+ \eta; z \in U, \eta \in C\}$$
for some open neighborhood $\omega$ of $p$ in $M$ and some
truncated open cone $C$ in $N_{p}(M)$, i.e. the intersection of a
open cone with a ball centered at 0.
Let $K$ be a proper closed subset of $M$, $K$ is said
{\it  $\cal W$-removable
at $p$} if there exists a wedge $\cal W$ of edge $M$ at $p$ such that
any CR function defined on $M\backslash K$ extends holomorphically to
$\cal W$.\par\noindent
The theory of CR removability theory has been first developed by the
deep work of Jor\"{\i}cke
\cite{joricke1,joricke2,joricke3}.
In the hypersurface case, she proved
the removability of proper closed subsets of 2 codimensional submanifolds of
$M$. Then Chirka-Stout \cite{chirkastout} proved the removability
of closed subsets of two codimensional null Hausdorff measure.
The CR removability results in the greater codimensional case was
obtained by Merker\cite{merker} and Joricke \cite{joricke3}. In
\cite{merker2}, Merker and Porten proved the removability of proper closed
subset of two codimensional submanifolds of $M$ and of closed subset
of $M$ of finite three codimensional Hausdorff measure.
All those removability theorems have been given under minimality
assumption on $M$.
Our main theorem gives a removability result for closed subsets of
two codimensional Hausdorff null measure under convexity
assumption given by the Levi form.
\medskip
\par\noindent
{\bf Main theorem}
{\it
Suppose $M$ is a generic, CR submanifold of an open set $U$ of $\CC^n$
 of class $C^{4}$
with $\dim_\RR M = 2n-d$ $(1 \leq d \leq
n-1)$. Let $p_0$ a point of $M$ such that the convex hull of the
image of the Levi form has nonempty interior.
Then their exists a wedge $\cal W$ of edge $M$ at $p_0$ such that if $K$ is
a closed subset of $M$ of null Hausdorff
${\cal H}^{2n-d-2}$ measure, every CR function $f$ on $M
\backslash K$ extends holomorphically to $\cal W$.
}\par
\medskip
\noindent
In the case $K$ is empty, this is the wedge of the edge extension
theorem of Airapetyan-Henkin \cite{ayrapetianhenkin} and
Boggess-Polking\cite{Boggess2}.\par\noindent 
Let $f$ be a CR-meromorphic map defined on a minimal generic
manifold $M$
and $K$ be the indeterminacy set of $f$.
If $x \not \in K$, as explained in \cite{fredy}, $f$ extends meromorphically
to a wedge of edge $M$ at $x$. Then, by the uniqueness of the extension,
all the meromorphic extensions of $f$ coincides. Thus, deforming $M$
outside $K$ in the wedge where the extension arise, we can assume that
  $f$ is the restriction of a meromorphic map defined in a
  neighborhood $U$ of $M \backslash K$.
By Oka-Levi theorem, the envelope of meromorphy of an open subset
of $\CC^n$ is the same than its envelope of holomorphy. So it suffice
to prove that the envelope of holomorphy of $U$
contains a wedge $\cal W$ of edge $M$ to prove that $f$ extend
meromorphically
into $\cal W$. As remarked in \cite{fredy}, $K$ is included (and of
empty inside) in a scarred two codimensional submanifold of $M$.
So applying our theorem and a removability theorem of \cite{merker2}
  for proper closed subsets of two codimensional submanifolds of $M$ we
obtain:
\begin{corol}
Suppose $M$ is a generic, $C^4$,  CR submanifold
of an open subset
$U$ of $\CC^n$ with $\dim_\RR M = 2n-d$ $(1 \leq d \leq
n-1)$.  Let $p_0$ a point of $M$ such that the convex hull of the
image of the Levi form has nonempty interior.
Then there exists a wedge $\cal{W}$ of edge $M$ at
$p_0$ such that every CR-meromorphic function defined on $M$ extends
meromorphically to $\cal{W}$.
\end{corol}
The authors are very greatful to Gennadi Henkin for proposing
them this problem. Also, the second author would like to thank J. Merker and
E. Porten for the
very instructive discussions he had with them
on the technics of analytic discs.
\section{Bishop's equation}
{\it 2.1 Normal form for CR submanifolds - }
\begin{lemme}
Suppose $M$ is a generic, CR submanifold of $\CC^n$ of class $C^k$
($k\geq 2$) with $\dim_\RR M= 2n-d$ ($1 \leq d \leq n$). Suppose that
$p_0$ is a point in $M$. There is a neighborhood $U$ of $p_0$ in
$\CC^n$ such that for any $p \in U \cap M$ there exists a biholomorphism
$\Phi=\Phi_p: U \rightarrow \Phi\{U\} \subset
\CC^n$; and a function $h=h_p: \RR^d \times \CC^{n-d} \rightarrow \RR^d$
of class $C^k$ with $h(0)=0$ ($\Phi$ and $h$ depending
in $C^{k-1}$ fashion of $p \in U \cap M$) such that
$$ \Phi \{ M \cap U\}= \{(x+iy,w)\in \Phi \{U\} \subset \CC^d \times
\CC^{n-d}; y=h(x,w)\}$$
Furthermore
$${\partial ^{|\alpha|+|\beta|}h(0)\over \partial
  x^\alpha\partial w^\beta}={\partial ^{|\alpha|+|\beta|}h(0)\over \partial
  x^\alpha\partial\bar w^\beta}=0$$ for all $\alpha \in \NN^d, \beta
  \in \NN^{n-d}$ such that
  $0 \leq|\alpha|+|\beta| \leq 2$.
\end{lemme}
This lemma can be found in \cite{Boggess}, we just remark that the
given construction of $\Phi$ and $h$ depends in $C^{k-1}$ fashion of
$p \in U \cap M$.
We will call the {\it quadric associated to $h$} the quadric
$q:={\partial^2h \over \partial w \bar \partial \bar w}(0)$.\par
\medskip
\noindent {\it 2.2 Bishop's equation - } Let $M$ be a generic CR submanifold
of
a neighborhood $U=U_1 \times U_2$ of the origin of $\CC^n$
defined by the equation
$$ M = \{ (z=x+iy,w): (x,w) \in U_1, y \in  U_2; y=h(x,w))\}$$
where $U_1$ is an open neighborhood of the origin of $\RR^d\times
\CC^{n-d}$, $U_2$ is an open neighborhood of the origin of ${\RR^d}$,
$h: U_1 \rightarrow U_2$ is of class $C^k$
and $h(0)=0, Dh(0)=0$. Let note $D$ the open unit disc of $\CC$,
$\bar D$ its closure and $S^1$ the unit circle of $\CC$.
Given an analytic disc $W: \bar D \rightarrow
\CC^{n-d}$, we wish to find an analytic disc $G: \bar D \rightarrow
\CC^d$ so that the boundary of the disc $A=(G,W): \bar D \rightarrow
\CC^n$ is contained in $M$. This mean that $G$ must satisfy
$$ \Im G(\zeta)=h(\Re G(\zeta),W(\zeta)) \mbox{ for } \zeta \in S^1. $$
This equation involves both $u= \Re G$ and $v = \Im G$. The above
equation will be easier to solve by eliminating $v$. To
do this, we use the Hilbert transform which is defined as follows.
If $u: S^1 \rightarrow \RR^d$ is a
continuous function, then $u$ extends to a unique harmonic function on the
unit disc $D$. This harmonic function has a unique harmonic conjugate
in $\bar D$ which vanishes at the origin. The Hilbert
transform of $u$ (denoted $Tu:S^1 \rightarrow \RR^d$) is defined to be
$v|_{S^1}+c$ where $c=-v(\zeta=0)$.
The function $-iG=v-iu$ is also analytic so $T(v|_{S^1})=-u+x$ where
$x=u(\zeta=0)$. Conversely, if $u,v: S^1 \rightarrow \RR^d$ are
continuous functions with $u=-Tv+x$, then $u+iv: S^1 \rightarrow
\CC^d$ is the boundary values of a unique analytic disc $G: \bar D
\rightarrow \CC^d$ with $\Re G(\zeta=0)=x$.
Suppose $u+iv=G: \bar D \rightarrow \CC^d$ is an analytic disc with
$v(e^{i\varphi})=h(u(e^{i\varphi}),W(e^{i\varphi}))$ for $0 \leq
\varphi \leq 2\pi$. We apply $-T$ to both sides of this equation and
obtain
$$u(e^{i\varphi})=-T(h(u,W))(e^{i\varphi})+x \mbox{ for } 0 \leq \varphi
\leq
2\pi$$
where $x \in \RR^d$ is the value of $u$ at $\zeta=0$. The above
equation will be referred to as Bishop's equation.
Conversely, suppose the analytic disc $W: \bar D \rightarrow
\CC^{n-d}$ and the vector $x \in \RR^d$ are given, and suppose $u:S^1
\rightarrow \RR^d$ is a solution to Bishop's equation. From the above
discussion, the function
$$ \varphi \rightarrow
u(e^{i\varphi})+ih(u(e^{i\varphi}),W(e^{i\varphi}))$$
is the boundary values of a unique analytic disc $G: \bar D
\rightarrow \CC^d$. Since $\Re G(e^{i\varphi})=u(e^{i\varphi})$, the
boundary of the analytic disc $A=(G,W):\bar D \rightarrow \CC^n$ is
contained in $M$. Furthermore, $\Re G(\zeta=0)=x$.
\par \noindent
The Hilbert Transform is a smooth linear map from the space
$C^\alpha(S^1,\RR^d)$ of $\alpha$-Hölderian maps from $S^1$
into $\RR^d$ to itself (with $0 < \alpha <1)$
(a proof can be found in \cite{Boggess}).\par

\section{The convex quadric case}
A quadric $M$ is a submanifold given by the equation
$$ M= \{(z=x+iy,w) \in \CC^d \times \CC^{n-d}; y=q(w,\bar w)\}$$
where $q:\CC^{n-d} \times \CC^{n-d} \rightarrow \CC^d$ is a quadric
form.
Let us start with a given analytic disc $W: \bar D \rightarrow
\CC^{n-d}$. The analytic disc $W: \bar D \rightarrow \CC^{n-d}$ is
given by a convergent power series
$$W_t(\zeta)=\sum_{j=0}^\infty t_j a_j \zeta^j,
t_j \in \RR, a_j \in \CC^{n-d}, \zeta \in \bar D.$$
In our application, all but a finite number of the parameters
$\{a_0,a_1,...\}$ and $\{t_0,t_1,...\}$ will vanish.
In order for the set $\{A(\zeta)=(G(\zeta),W(\zeta)); \zeta \in S^1 \}$ to
be contained in $M$, the analytic disc $G: \bar D \rightarrow \CC^d$
must satisfy
$$\Im G(\zeta)=q(W(\zeta),\overline{W(\zeta)}) \mbox{ for } \zeta \in S^1.$$
Using the linearity and symmetry of $q$, we find
$$G(\zeta)=x+i\sum_{j=0}^\infty t_j^2 q(a_j,\bar a_j) +2i \sum_{0 \leq k <
  j} t_jt_k q(a_j,\bar a_k)\zeta^{j-k}$$
(with $x \in \RR^d$) is, given the condition $\Re G(0)=x$,
 the unique solution of the equation $\Im
  G(\zeta)=q(W(\zeta),\overline{W(\zeta)})$ for $\zeta \in S^1$.
Then we have
$$A(\zeta) \in M \mbox{ for } \zeta \in S^1 $$
$$A(\zeta=0)=\left(x+i\sum_{j=0}^\infty t_j^2 q(a_j,\bar a_j),t_0a_0
\right)$$
We want to verify that when $t$ moves, the boundaries of the obtained
discs  give us a submersion of the space of parameters.
Take $t_0=1$.
It suffice to verify that the matrix $\cal M$ of the derivatives
of $(W,\bar W, \Re G,v(0))$ (we
recall that  $v(0)=\sum_{j=0}^\infty t_j^2 q(a_j,\bar a_j) $) in $t_i$
is of maximal rank.
We have for $\zeta \in S^1$
$${\partial \over \partial t_j}v(0)= 2t_jq(a_j,\bar a_j)$$
and
\begin{eqnarray*}
\displaystyle{\partial \over \partial t_j}\Re G(\zeta) & = &
\Re\left( 2i \sum_{1 \leq k < j} t_k q(a_j,\bar a_k)\zeta^{j-k}
+ 2i\sum_{k > j} t_kq(a_k,\bar a_j)\zeta^{k-j}  \right) \\
& = & i\left(
\sum_{k<j} t_kq(a_j,\bar a_k)\zeta^{j-k}
+\sum_{k>j} t_kq(a_k,\bar a_j)\zeta^{k-j}-\right.\\
& &\hspace{2cm}
\left.- \overline{\sum_{k<j} t_kq(a_j,\bar a_k)\zeta^{j-k}}
-\overline{\sum_{k>j} t_kq(a_k,\bar a_j)\zeta^{k-j}}
\right)\\
& = & i\left(
\sum_{k<j} t_kq(a_j,\bar a_k)\zeta^{j-k}
+\sum_{k>j} t_kq(a_k,\bar a_j)\zeta^{k-j}-\right.\\
& & \hspace{2cm}\left. -\sum_{k<j} t_kq(a_k,\bar a_j)\zeta^{k-j}
-\sum_{k>j} t_kq(a_j,\bar a_k)\zeta^{j-k}
\right)
\end{eqnarray*}
We remark that the vectors $\{ i(\sum_{k=1}^\infty t_kq(a_j,\bar
  a_k)\zeta^{j-k}-\sum_{k=1}^\infty t_kq(a_k,\bar
a_j)\zeta^{k-j})\}_{j\geq 1}$ are linear combinations of the vectors
  $\{a_j\zeta^j\}_{j \geq 1}$
and of the vectors $\{\bar a_j \zeta^{-j}\}_{j \geq 1}$
(we recall that the $a_j$ are vectors and that $q$ takes vectorial
values).
Then by adding this vectors to the derivative of
$\Re G(\zeta)$ we does not change the rank of the matrix of the
derivatives in $t_j$.
So the rank of the  matrix $\cal{M}$ of the derivatives of $(W, \bar
W, \Re G,v(0))$ in
$t_j$, with $\zeta \in S^1$
$$ {\cal M}(a,\zeta,t)= \left(
\begin{array} {cccc}
a_1\zeta^1 & \ldots &a_j\zeta^j & \ldots \\
 \bar a_1\zeta^{-1} &\ldots & \bar a_j\zeta^{-j} & \ldots \\
\displaystyle {\partial \over \partial t_1}\Re G(\zeta) & \ldots &
 \displaystyle {\partial \over \partial
  t_j}\Re G(\zeta) & \ldots \\
 2t_1q(a_1,\bar a_1)&\ldots & 2t_jq(a_j,\bar a_j) & \ldots
\end{array}
\right)$$
is of the same rank than the matrix
$${\cal  M}'(a,\zeta,t)=\left(
\begin{array} {cccc}
a_1\zeta^1 & \ldots &a_j\zeta^j & \ldots \\
 \bar a_1\zeta^{-1} &\ldots & \bar a_j\zeta^{-j} & \ldots \\
 P_1(a,\zeta,t) & \ldots &  P_j(a,\zeta,t) & \ldots \\
 t_1q(a_1,\bar a_1)&\ldots & t_jq(a_j,\bar a_j) & \ldots
\end{array}
\right)$$
where $$P_j(a,\zeta,t)=2\sum_{k<j}
\left(
q(a_j,\bar a_k)\zeta^{j-k} - q(a_k,\bar a_j)\zeta^{k-j} \right)$$
is a polynomial in $\zeta$ and $\zeta^{-1}$ of degree $j$ and depends
only o
n $(t_1,...,t_{j-1})$.
Let note ${\cal M}'_N$ the matrix ${\cal M}'$ truncated at the $N^{th}$
column.
\begin{lemme}
For any fixed $\zeta \in S^1$ there exists
$N \in \NN$, $b=(b_j)_{1 \leq j \leq N}$ and
$s=(s_j)_{1 \leq j \leq N}$ such that the rank of ${\cal M}'_N(b,\zeta,s)$
is maximal.
\end{lemme}
\begin{demo}
First, let suppose $\zeta=1$. Let suppose, that for all  $N \in \NN$,
$b=(b_j)_{1 \leq j \leq N}$ and
$s=(s_j)_{1 \leq j \leq N}$ the rank of  ${\cal M}'_N(b,1,s)$ is
not maximal. Let suppose $N,b, s$ are chosen such that the rank
$r$ of
${\cal M}'_N(b,1,s)$ is the greater rank reached by those matrixes.
As the rank is never maximal we have $r < 2n$.
 Then there exists a non null $(1,2n)$ matrix
$(\alpha,\beta,\gamma,\delta) \in \CC^{n-d}\times \CC^{n-d} \times
\CC^d \times \CC^d$ such that $(\alpha,\beta,\gamma,\delta)
\times{\cal M}'_N(b,1,s)=0$.
So for any families  $b=(b_j)_{N \leq j \leq N'}$ and
$s=(s_j)_{N \leq j \leq N'}$, we have
$(\alpha,\beta,\gamma,\delta) \times{\cal M}'_{N'}(b,1,s)=0$
because $ {\cal M}'_N(b,1,s)$ is a sub-matrix of ${\cal
  M}'_{N'}(b,1,s)$.
As the $2n-d$ first lines of the matrix are independent of
$s_{N'}$, we have $\delta \times (s_{N'}q(b_{N'}, \bar b_{N'}))=0$ for
all $N'>N$, $s_{N'} \in \RR$ and $b_{N'} \in \CC^{n-d}$.
From the hypothesis that the convex hull of the image of $q$ is an
open cone of $\RR^d$, the previous equality implies that $\delta=0$.
As the $2n-2d$ first lines of the matrix are independent of
$s_{N'-1}$, we have $\gamma \times s_{N'-1}(q(b_{N'},\bar b_{N'-1})
- q(b_{N'-1},\bar b_{N'}))=0$ for
all $N'>N+1$, $s_{N'-1} \in \RR$ and $b_{N'-1},b_{N'} \in \CC^{n-d}$.
Choosing for example $b_{N'}=ib_{N'-1}$
we find $q(b_{N'},\bar b_{N'-1})
- q(b_{N'-1},\bar b_{N'})=2iq(b_{N'-1},\bar b_{N'-1})$.
So, from the hypothesis that the convex hull of the image of $q$ is an
open cone of $\RR^d$  the equality implies that $\gamma=0$.
Then, it is obvious that $\alpha=\beta=0$ and we have
the contradiction.\par\noindent
Let choose $k \in \NN$, $b^*=(b_j^*)_{1 \leq j \leq k}$ and
$s^*=(s_j^*)_{1 \leq j \leq k}$ such that the rank of ${\cal
  M}'_k(b^*,1,s^*)=2n$, i.e. is maximal.
Let $\zeta$ fixed, there exists $r \in \NN^+$ such that
$|\zeta^{mr}-1|<\epsilon$ for all $1 \leq m \leq k$ for $\epsilon >0$
small enough.\par\noindent
For all $1 \leq m \leq k$, let choose $b_{mr}=b_m^*$ and
$s_{mr}=s^*_m$ and $b_j=s_j=0$ for all $ j \not \in \{r,2r,...,kr\}$.
Then, by continuity, for $\epsilon$ small enough we have
$$\mbox{rank} {\cal M}'_{kr}(b,\zeta,s)= \mbox{rank} {\cal
M}'_k(b^*,1,s^*)=2n.$$
\end{demo}
\begin{prop}
There exists  $N \in \NN$, $a=(a_j)_{j \leq N}$ and an open cone
$\Omega \subset \RR^N$ such the rank of ${\cal M}(a,t,\zeta)$
is maximal for any $t \in \Omega$ and  $\zeta \in S^1$.
\end{prop}
\begin{demo}
We just remark that $\mbox{ rank } {\cal M}'_N(a,\zeta,t)=
\mbox{ rank } {\cal M}'_N(a,\zeta,\lambda t)$ for all $\lambda \in
\RR\backslash\{0\}$. So it is enough to prove our proposition for
$\Omega$ any non-empty open set in $\RR^n$ and by continuity it
suffice to find $N, b$ and $t$ fixed such that for any $\zeta \in
S^1$, $\mbox{ rank } {\cal M}'_N(a,\zeta,t)=2n$.\par \noindent
For all $\zeta \in S^1$, let choose $N_\zeta, \{b_j(\zeta)\}_{1 \leq j \leq
  N_\zeta}$ and $  \{s_j(\zeta)\}_{1 \leq j \leq N_\zeta}$ verifying the
previous lemma. By continuity, there exists a neighborhood $V_\zeta$
of $\zeta$ in $S^1$ such that
$\mbox { rank } {\cal M}'_{N_\zeta}(b_j(\zeta), \zeta',s(\zeta))=2n$ for all
$\zeta' \in V_\zeta$.
Let $\{\zeta_1,...,\zeta_k\} \subset S^1$ a finite subset such that
$S^1 \subset \bigcup_{m=1}^k V'_{\zeta_m}$ where $V'_{\zeta_m}$ is a
relatively compact open subset of $V_{\zeta_m}$.
Let $N=N_{\zeta_1}+...+N_{\zeta_k}$,
$a_j=b_{j-(N_{\zeta_1}+...+N_{\zeta_m})}(\zeta_m)$,
$t_j=\theta_m s_{j-(N_{\zeta_1}+...+ N_{\zeta_m})}(\zeta_m)$
for all $N_{\zeta_1}+...+N_{\zeta_m} < j \leq
N_{\zeta_1}+...+N_{\zeta_{m+1}}$, $0 < \theta_m << \theta_{m+1}$ and $0 \leq
m \leq k-1$.
There is still to show that $\mbox{rank} {\cal M}'_N(a,\zeta,t)=2n$
for all $\zeta \in S^1$. Let $\zeta \in V'_{\zeta_m}$, we have
$\mbox{rank } {\cal M}'_N(a,\zeta,t) \geq \mbox{ rank }{\cal
  M}'_{N_1+...+N_m}(a,\zeta,t)$.
The rank of the second matrix is greater or equal to the one of the
sub-matrix formed by the $N_m$ last columns. As the $\theta_m$ are
chosen such that $\theta_{m-1} <<\theta_m$, the rank of this
sub-matrix is equal to the rank of ${\cal
  M}'_{N_m}(b(\zeta_m),\zeta,s(\zeta_m))$, i.e. is maximal.
\end{demo}
\par
\medskip
\noindent We note $A(p,a,t)=(G(p,a,t),W(p,a,t))$ the attached discs
given by the previous proposition where $p=(x,a_0) \in \RR^d \times
\CC^{n-d}$, $a=(a_1,...,a_N)$ and $t=(t_1,...,t_N)$.
From \cite{Boggess} the image, for all $t \in \Omega$, of
$A(p,a,t)(0)$ contains a wedge ${\cal W}=\omega+{\cal C}$ of edge $M$ at
$0$.
Let $\Omega_z$ be the set of $t \in \Omega$ such that $A(p,a,t)(0)=z$
 where $z=p+\eta$ with $p \not \in K$, $\eta \in {\cal C}$
and  $K$ a closed subset of $M$ verifying the hypothesis of the theorem 1.
From our proposition, the map $(u,W)(p,a,.)(.): \Omega_z \times S^1
\rightarrow \RR^d \times \CC^{n-d}$ is a submersion. So, their exists
discs, attached to $M$, passing through $z$, isotopic to the point $p$ such
that their boundary is in $M \backslash K$. Here, we say that a disc
$A$ is isotopic to a point $p \in M \backslash K$,
if there exists a continuous family of discs attached to $M \backslash
K$ and of class $C^2$
that contains $A$ and $p$.
This implies the theorem 1 in the quadric case. The details can be
find at the end of the next section.

\section{The convex case}
Let $p_0 \in M$, $U$, $\Phi=\Phi_p$ and $h=h_p$ verifying  lemma 1.
Let $U^p:=\Phi_p(U)$, $M^p:=\Phi_p(M \cap U)$, $K^p=\Phi_p(K)$
, $\Pi_1: \CC^n \rightarrow
\RR^d \times \CC^{n-d}$ be the projection $\Pi_1(u+iv,w)=(u,w)$ and
$\Pi_{2}:\CC^n \rightarrow \RR^d$ be the projection
$\Pi_{2}(u+iv,w)=v$. Let $\widetilde \Pi: \RR^d
\backslash 0 \rightarrow S^{d-1}$ be the usual projection into the
unit sphere $S^{d-1}$ of $\RR^d$ and $\Pi':=\Pi_1|_M$.
Let $N \in \NN$, $a_j \in \CC^{n-d}$ and $t_j \in \RR^N$ fixed. We consider
the
analytic disc $W(t,.): \bar D \rightarrow \CC^{n-d}$ with
$W(t,\zeta)=\sum_1^N t_ja_j\zeta^j$. Then if $|t_ja_j|$ is small enough,
Bishop's equation (for $M^p$ with the condition $u(0)=W(0)=0$) admits a
unique solution
depending in a $C^4$ fashion of $t_j$ \cite{Boggess}.
Let $${\cal M}:=\displaystyle \left({\partial W(t,\zeta) \over
    \partial t_j },
{\partial \bar W(t,\zeta) \over \partial t_j},
{\partial u(t,\zeta) \over \partial t_j},
{\partial v(t,0) \over \partial t_j}\right)_{1 \leq j \leq N}$$
\begin{prop}
Let suppose that $\Gamma_{p_0}$ has non empty interior. Then
there exists $\omega$ a neighborhood of $p_0$ in $M$, $N \in \NN$,
$a_j \in \CC^{n-d}$, $\Omega \subset \RR^N$ a truncated open cone such
that for every $p \in \omega$, $t \in \Omega$, $\zeta \in S^1$ the
rank of the matrix ${\cal M}$ is maximal.
\end{prop}
\begin{demo}
We call $q$ (resp. $q_0$) the quadric form associated to $h=h_p$
(resp. $h_{p_0}$) defined in the section 2.1.
Let define ${\cal O}^\alpha(\bar D,\CC^{n-d})=C^\alpha(\bar D,\CC^{n-d})
\cap {\cal O}(D,\CC^{n-d})$ the space of holomorphic maps from $D$
into $\CC^{n-d}$ which are $\alpha$-hölderian up to the boundary
and $$\displaystyle  {\cal B}(\bar B^{2n-d},\RR^d)=
C^k(\bar B^{2n-d},\RR^d)\cap \{ {\partial ^{|\alpha|+|\beta|}g \over
  \partial ^\alpha x \partial ^\beta w}={\partial ^{|\alpha|+|\beta|}g \over
  \partial ^\alpha x \partial ^\beta \bar w}=0, 0 \leq |\alpha|+|\beta| \leq
2 \}$$ where $\bar B^{2n-d}$ is the closed unit ball of $\RR^{2n-d}$.
Let consider the map $${\cal F}: C^\alpha(\bar D, \RR^d)\times {\cal
  O}^\alpha(\bar D, \CC^{n-d})\times {\cal B}(\bar B^{2n-d},\RR^d)
\rightarrow C^\alpha(\bar D, \RR^d)$$
$$ {\cal F}(u,W,g):=u+T(g(u,W))$$
For  $0<\alpha<1$, ${\cal F}$ is of  class $C^1$ and $$\displaystyle
{\partial {\cal F} \over \partial u}(0,0,q_0)={\rm Id}$$
From the implicit functions theorem, in a neighborhood of $(0,0,q_0)$,
their exists a unique solution $u(W,g)$ of the equation ${\cal
  F}(u,W,g)=0$ depending in a $C^3$
fashion of $W$ and $g$.
Let ${\cal M}_g$ be the matrix defined as above for the manifold
$\{y=g(x,w)\}$
 and $\bar \Theta$ be a compact of non empty inside $\Theta$  of
$\Omega_0 \cap S^{N-1}$ where  $S^{N-1}$ is the unit sphere of $\RR^N$
and $\Omega_0$ and  $N$ are defined in the proposition 1 for $q=q_0$.
Then for
$|g-q_0|< \epsilon$ with $\epsilon$ small enough and $t \in \bar \Theta$,
the rank of ${\cal M}_g$ is still maximal.
Let $\omega$ be a neighborhood of $p_0$ in $M$ small enough such that for
$p \in \omega$ we have $|q-q_0|<\epsilon/2$. Let $0 < \lambda <1$
and set $h_\lambda(x,w):={1\over \lambda^2}h(\lambda x,\lambda
w)$. From lemma 1,  their exists $\delta > 0$ such that
$|h_\lambda-q|< \epsilon/2$ for  $0 < \lambda
<\delta$. So, for all $t \in \bar \Theta$, the rank of
${\cal M}_{h_\lambda}$ is maximal for the submanifold
$\{y=h_\lambda(x,w)\}$.
So, for all $t \in \lambda \Theta$, the rank of $\cal M$
is maximal for ${M}^p$. Then it suffice to take $\Omega=\{t \in \lambda
\Theta ; 0 < \lambda < \delta\}$.
\end{demo}\par
\medskip
\noindent From lemma 1, $\{ \Phi_p^{-1} \}_{p \in \omega}$ is a family of
holomorphic maps depending of $C^{k-1}$ fashion of $p \in \omega$. So,
the coefficients of the Taylor expansion of $\Phi_p^{-1}$ at $0$ also
depends of $C^{k-1}$ fashion of $p \in \omega$. Without loss of
generality, we can assume that $p_0=0$ and that
$T_M(p_0)=\{\Im u=0\}$.
Let $\Omega \subset
S^{d-1}$ be a small neighborhood of a point $\theta \in S^{n-d}$.
The map
$$F=(F_1,F_2,F_3): ]-\epsilon,\epsilon[ \times \Omega \times \omega
\rightarrow
\RR \times S^{d-1}\times M$$
defined by
$$F_3 (\lambda, \theta,p)={\Pi'}^{-1}\circ\Pi_1\circ
\Phi_p^{-1}(i\lambda\theta,0)$$
$$F_1(\lambda, \theta,p)=
\left\{{\begin{array} {ll}
|\Pi_2(\Phi_p^{-1}(i\lambda\theta,0)-F_3(\lambda,\theta,p))|
& \mbox{ if } \lambda > 0 \\
-|\Pi_2(\Phi_p^{-1}(i\lambda\theta,0)-F_3(\lambda,\theta,p))|
& \mbox{ if } \lambda < 0
\end{array}}
\right.$$
$$F_2 (\lambda, \theta,p)=
\left\{{\begin{array} {ll}
\widetilde{\Pi}(
\Phi_p^{-1}(i\lambda\theta,0)-F_3(\lambda,\theta,p))
& \mbox{ if } \lambda > 0 \\
-\widetilde{\Pi}(
\Phi_p^{-1}(i\lambda\theta,0)-F_3(\lambda,\theta,p))
& \mbox{ if } \lambda < 0
\end{array}}
\right.$$
for $\epsilon > 0$ small enough and $\lambda \neq 0$, extends
at $\lambda =0$ as a $C^{1}$ map. Indeed, we have $\Phi_p^{-1}=L+r$
with $L$ a linear map and $r \in {\cal O}(\lambda^2)$ when $\lambda$
tends to $0$. As for
$\Phi_p^{-1}$ replaced by $L$,
$F$ extends at $\lambda=0$ as a $C^1$ map, this is
also true for $\Phi_p^{-1}$.
\begin{lemme}
For every truncated open cone $C \subset \RR^d$
$$\bigcup_{p \in \omega} \Phi_p^{-1}(\widetilde C)$$
contains a wedge $\cal W$ of edge $M$ at $p_0$ where $\widetilde C
:=\{z=(i\eta,v)\in \CC^d \times \CC^{n-d}\}$. Moreover, this union
form a sheeting of class $C^{k-1}$ of $\cal W$.
\end{lemme}
\begin{demo}
Let $\theta \in S^{d-1}$ fixed and  $\Omega \subset S^{d-1}$ be a
small enough neighborhood of $\theta$ such that $\delta \Omega \subset C$
for
$\delta > 0$ small enough.
So, in a neighborhood of the point $(0,
\theta,p_0)$, $F$ is homeomorphic to its image for all $\theta \in
S^{d-1}$. This is obvious because
${\Pi'}^{-1}\circ\Pi_1\circ \Phi_p^{-1}(0,.)$ is the
identity map of $\omega$ and $\widetilde
F(.,.,p_0):=(F_1(.,.,p_0),F_2(.,.,p_0))$  is a homeomorphism
because $\Phi_{p_0}$ is a conformal map (as $\Phi_{p_0}$ is a
biholomorphism).
 As the image of $F(0,.,.)$ is in
$\{0\}\times S^{d-1}\times M$,
$F(]0,\epsilon[\times\Omega\times\omega)$ contains
$V^+:=V \cap \{\RR^+\times S^{d-1}\times M\}$ where $V$ is a
neighborhood of $F(0,\theta,p_0)$, we have that  $\bigcup_{p \in \omega}
\Phi_p^{-1}(\widetilde C)$ contains the set
$${\cal W}':=\{p+\eta\mbox{  with } \eta=\lambda \theta \mbox{ and }
(\lambda,\theta,p) \in V^+\}$$ which contains a wedge $\cal W$ of edge $M$.
\end{demo}
\begin{corol}
There exists a wedge $\cal W$ of edge $M$ at $p_0$ such that for all
$K \subset M$ (not necessary closed) of null
${\cal H}^{2n-d-1}$ measure and for every $z \in {\cal W}$ their exists
an analytic disc attached to $M$ of class $C^k$ passing through $z$ and
which border does not meet $K$.
\end{corol}
\begin{demo}
From lemma 3, it suffice to prove that for a fixed point $p \in
\omega$ and for all $z \in \{A(t,0): t \in
\Omega\}$, there exists an analytic disc attached to $M^p$, of class
$C^3$, passing through $z$ and whose border does not meet $K^p$.
Where $A(t,\zeta):=(G(t,\zeta),W(t,\zeta))$ is the analytic disc
attached to $M^p$ obtained by the Bishop's equation with
$({\rm Re} G, W)=(0,0)$. Indeed, from \cite{Boggess} $\{A(t,0); t \in
\Omega\}$
contains an open cone of $\RR^d$.\par\noindent
Let $\Theta$ and $\delta$  be as defined here above . Let define  the map
$$\Psi: ]0,\delta[\times \Theta \times S^1 \rightarrow \RR^d\times
\RR^d\times \CC^{n-d}$$
$$
\Psi(\lambda,\theta,\zeta)=(v(\lambda\theta,0),u(\lambda\theta,\zeta),
W(\lambda\theta,\zeta)).$$
where $G=u+iv$.
From the proposition 2, for all $\eta \in \{v(t,0); t \in \Omega\}$
$$E_\eta:=\{(\lambda,\theta) \in ]0,\delta[\times \Theta \mbox{ such that }
v(\lambda\theta, 0)=\eta\}$$
is a submanifold of codimension $d$ of $]0,\delta[\times\Theta$, i.e. of
dimension $N-d$ and $$\Psi_\eta:E_\eta \times S^1 \rightarrow \RR^d \times
\CC^{n-d}$$
$$\Psi_\eta(\lambda,\theta,\zeta)=(u(\lambda \theta,\zeta),W(\lambda
\theta,\zeta))$$ is a
submersion. We remark that the map $\Psi_\eta(\lambda,\theta,.)$ is
the projection of the border of the attached disc on $\RR^d \times
\CC^{n-d}$.
It follows that  ${\cal H}^{N-d}(\Psi_\eta^{-1}(\Pi_{1}(K^p)))=0$
and (from \cite{federer} 2.10.25) that for ${\cal H}^{N-d}$-almost all
$(\lambda,\theta) \in E_\eta$,
$\Psi_\eta(\lambda,\theta,S^1) \cap \Pi_{1}(K^p)= \emptyset$.
So, for such a $(\lambda,\theta)$, the attached disc does not meet $K$.
\end{demo}
\begin{corol} Let $K$ be a compact subset of $M$ of null
$2n-d-2$ dimensional Hausdorff measure. Then there exists
an open truncated cone $C \subset \RR^d$ such that
for all $z \in \Phi_p^{-1}(\{0\}\times C \times \{0\})$ with
$ p \not \in K$, there exists a  analytic disc attached to
$M\backslash K$ passing through $z$ and isotopic to $p$.
\end{corol}
\begin{demo}
With the same notations of the proof of the previous corollary,
from lemma 1 we have that ${\cal H}^{N-1}(\bigcup_{\eta \in \RR^d}
\Psi_\eta^{-1}(\Pi_1(K^p))=0$. So, $\bigcup _{\eta\in \RR^d}
\Psi_\eta^{-1}(\Pi_1(K^p))$ does not disconnect $]0,\delta[\times \Theta$.
So, for all $p \in \omega \backslash K$ and $z \in \{A(t,0), t \in
\Omega\}$ there exists an analytic disc attached to $M^p\backslash
K^p$ passing through $z$ and isotopic to $0$ and this
implies the corollary because $\Phi_p$ is a local biholomorphism.
\end{demo}\par
\medskip\noindent
{\it Proof of theorem 1. }
 From corollary 3 and the continuity principle
\cite {merker}, $f$ extends holomorphically in a neighborhood of each
point of $\bigcup_{p \in \omega \backslash K} \Phi^{-1}_p(\widetilde C)$.
Moreover, this extension is univalued. Indeed, the value of the
extension is given by the Cauchy formula and is locally constant in
the space of parameters $E_\eta$ of the attached discs. As $K$ is of null
Hausdorff codimension 2, the set of parameters of which attached discs
intersect $K$ does not disconnect $E_\eta$. 
For $\lambda$ small enough, the situation is close to the
quadratic case, so there exists $\epsilon$ such that
$E_\eta^\lambda=E_\eta \cap \{|\lambda|<\epsilon\}$ is connected. So the
extension given by the Cauchy formula is univalent in all the
connected component of $E_\eta^\lambda$ in $E_\eta$. Doing so for all $\eta$ we
obtain the unicity of the extension.\par\noindent
Let $V=\bigcup_{p \in \omega}\Phi_p^{-1}(\widetilde C)$
and $S=\bigcup_{p \in K}\Phi_p^{-1}(\{0\}\times C \times \{0\})$ .
From lemma 3, ${\cal H}^{2n-2}(S)=0$
and $f$ extends holomorphically in $V \backslash S$. 
By Hartog's theorem, $f$ extends holomorphically in $V$. Lemma 3
proves that $V$ contains a wedge $\cal W$ of edge
$M$ at $p_0$ and independent of $f$.

\par
\bigskip
\par\noindent
T.C. Dinh \ \ \ F. Sarkis \par\noindent
{\it Institut de math\'ematiques\par\noindent
Universit\'e Pierre et Marie Curie, Tour 46-56 B. 507\par\noindent
4 Place Jussieu, 75005 Paris\par\noindent
e-mail:dinhtien@math.jussieu.fr\par\noindent
e-mail:sarkis@math.jussieu.fr}\par\noindent
\end{document}